\documentclass[12pt, a4paper]{article}
\usepackage{amsthm}
\usepackage{amsmath,amssymb,latexsym,amsfonts,mathrsfs}
\usepackage[dvipdfmx]{graphicx}

\usepackage{color}
\usepackage{mathrsfs}

\usepackage{fancyhdr}
\usepackage{braket}
\usepackage[margin=20truemm]{geometry}

\theoremstyle{definition}
\newtheorem{Thm}{{\bf Theorem}}[section]

\newtheorem{Lem}[Thm]{{\bf Lemma}}
\newtheorem{Prop}[Thm]{{\bf Proposition}}

\newtheorem{Def}[Thm]{{\bf Definition}}
\newtheorem{Rem}[Thm]{{\bf Remark}}

\newcommand{\R}{\mathbb R}

\numberwithin{equation}{section}
\title{Singularities of solutions to\\ nonlinear Schr\"{o}dinger equations }
\author{Fumihito Abe and Keiichi Kato}
\begin{document}
\maketitle
\begin{abstract}
The aim of this work is to study the wave front set of the solutions of the initial value problem for nonlinear Schr\"{o}dinger equations. We determine the $H^s$ wave front sets of solutions from the behavior at infinity of the initial data by using the characterization of the $H^s$ wave front set by wave packet transform.
\end{abstract}
\section{Introduction}
\indent
In this paper, we study the wave front sets of the solutions of the following initial value problem of the nonlinear Schr\"odinger equations,
\begin{align}
\begin{cases}
i\partial_tu+\frac{1}{2}\Delta u=\mathcal{N}[u],\qquad (t,x)\in\R\times\R^n,\\
\left.u\right|_{t=0}=u_0,\qquad x\in\R^n,
\end{cases}\label{1}
\end{align}
where $u:\R\times\R^n\rightarrow\mathbb{C}$ is unknown function, 
$\mathcal{N}[u]=u^p \bar{u}^q$ with $p, q\in \mathbb{N}$, 
$u_0\in H^{s,0}(\R^n)\cap H^{0,s}(\R^n)$ with $s>n/2$, 
$i=\sqrt{-1}$, 
$\partial_t=\frac{\partial}{\partial t}$ and $\Delta=\sum_{j=1}^n\frac{\partial^2}{\partial x_j^2}$. 
\par
Singularities of solutions to linear Schr\"odinger equations have been studied by many authors (\cite{Craig}, \cite{Hassell}, \cite{Nakamura1},\cite{Nakamura3},\cite{Nakamura2}, \cite{Okaji}, \cite{Okaji2}, \cite{Parenti}, \cite{Sakurai} ). 
\par
For nonlinear Sch\"odinger equations, J. Szeftel \cite{JS} has obtained the same result as in \cite{Craig} by W. Craig, T. Kappeler and W. Strauss. He has shown that the solutions is in $H^r$ in the direction $\xi_0$ for $t>0$ if the initial data is in $H^{s,0}\cap H^{0,s}$ with $s>n/2$ and decays in the direction $\xi_0$ of order $r$ with $s<r<2s -n/2$.
In this paper, we show that for given $t_0>0$ the solution $u(t_0,\cdot )$ to nonlinear Schr\"odinger equations is in $H^s$ in the direction $\xi_0$ if the initial
 data satisfies the condition \eqref{131814_19Jun24} and \eqref{1_July17} via wave packet transform which is weaker than the condition on the initial data in \cite{JS}. 
Our condition of the main result assures that the solution $u(t_0, \cdot )$ is in $H^r$ in the direction of $\xi_0$ and does not imply $H^r$ regularities of the solutions $u(t, \cdot )$ with $t\ne t_0$, while the result by Szeftel \cite{JS} says that $u(t,\cdot )$ in 
in $H^r$ in the direction of $\xi_0$ for all $t\ne 0$. 
\par
In order to state our results precisely, we prepare several notations and give the definition of wave packet transform and the $H^s$ wave front set. We write $\hat{f}(\xi)=\mathcal{F}[f](\xi)= (2\pi)^{-n/2}\int_{\R^n}f(x)e^{-ix\xi}dx$ for Fourier transform of $f$. Inverse Fourier transform ...For $x=(x_1,\ldots,x_n)\in\R^n$, we denote $\braket{x}=(1+|x|^2)^{1/2}$. For  $\varphi\in\mathcal{S}(\R^n)\backslash\{0\}$, we put $\varphi_{\lambda}(x)=\lambda^{nb/2}\varphi(\lambda^b x)$ with $0<b<1$ and $\lambda\geq1$, and write $\varphi^{(t)}_{\lambda}(x)=e^{it\Delta/2}\varphi_{\lambda}(x)$, where $e^{it\Delta/2}$ is the evolution operator of the free Schr\"{o}dinger equation. 
${\mathcal S}$ stands the set of all Schwartz functions (rapidly decreasing functions) and ${\mathcal S'}$ stands its dual. 
We denote the weighted Sobolev spaces $\left\{f\in L^2(\R^n) \mid (1-\Delta)^{s/2}[\braket{\cdot}^mf]\in L^2(\R^n)\right\}$ by $H^{s,m}$. 
For a subset $K$ in $\R^n$, $\overline{K}$ stands the closure of $K$. 
\begin{Def}[Wave packet transform]
Let $\varphi\in\mathcal{S}(\R^n)\backslash\{0\}$. For $f\in\mathcal{S'}(\R^n)$, the wave packet transform $W_{\varphi}f(x,\xi)$ of $f$ with window function $\varphi$ is defined as
\[
W_{\varphi}f(x,\xi)=\int_{\R^n}\overline{\varphi(y-x)}f(y)e^{-iy\xi}dy,\quad x,\xi\in\R^n.
\]
\end{Def}
\begin{Def}[$H^s$ wave front set]
For $f\in\mathcal{S'}(\R^n)$, $WF_{H^s}(f)$ is a subset in $\R^n\times(\R^n\backslash\{0\})$ is said to be $H^s$ wave front set if the following condition holds: \\
$(x_0,\xi_0)\notin WF_{H^s}(f)$ if and only if there exist a function $\chi\in C_{0}^{\infty}(\R^n)$ with $\chi(x_0)\neq 0$ and a conic neighborhood $\Gamma$ of $\xi_0$ in $\R^n\backslash\{0\}$ such that $\braket{\xi}^s |\widehat{\chi f}(\xi)|\in L^2(\Gamma)$.
\end{Def}
The following theorem is our main result.
\begin{Thm}
\label{2}
Let $\varphi (x)=e^{-x^2/2}$ and  a real number $s$ be larger than $n/2$ and take a positive constant $b\le 1/2$ sufficiently small which depends only on $s$, $n$, $p$ and $q$. Suppose that $u_0$ in $H^{s,0}(\R^n)\cap H^{0,s}(\R^n)$, 
that $u(t,x)$ is a solution to (\ref{1}) in $C([-T, T];H^{s,0}(\R^n)\cap H^{0,s}(\R^n))$ and that there exist positive numbers $t_0$ and $\delta$, and a neighborhood $V$ of $\xi_0$ such that
\begin{equation}
\label{131814_19Jun24}
 \int_1^{\infty}\lambda^{2r+n-1}\underset{z\in \R^n}{\mbox{sup}}\left\|W_{{\varphi}_{\lambda}^{(-t_0)}}u_0\left(x-t_0\lambda\xi, \lambda\xi\right)\right\|^2_{L^2_x(B(z,\delta)) \times L^2_{\xi}(V)}d\lambda<+\infty
\end{equation}
and
\begin{equation}
\label{1_July17}
 \int_1^{\infty}\lambda^{2r+n-1}\underset{z\in \R^n}{\mbox{sup}}\left\|W_{{\varphi}_{\lambda}^{(t_0)}}u_0\left(x+t_0\lambda\xi, \lambda\xi\right)\right\|^2_{L^2_x(B(z,\delta)) \times L^2_{\xi}(-V)}d\lambda<+\infty
\end{equation}
with some $\xi_0\in\R^n\backslash\{0\}$ and with some positive number $r$ satisfying $s<r<2s-n/2$. Then we have $(x,\xi_0)\notin WF_{H^r}[u(t_0)]$ for all $x\in \R^n$, where $\varphi_{\lambda}^{(t)}(x)=e^{it\Delta/2}\varphi_{\lambda}(x)$ with $\varphi_{\lambda}(x)=\lambda^{nb/2}\varphi(\lambda^bx)$
and $-V=\{\xi \in \R^n | -\xi \in V\}$. 
\end{Thm}
\begin{Rem}
In the case that $\mathcal{N}(u, {\bar{u}})=u^p$ with $p\in \mathbb{N}$, the assersion of Theorem \ref{2} holds without the condition \eqref{1_July17}. 
In the case that $\mathcal{N}(u, {\bar{u}})=\overline{u}^q$, the assersion holds without the condition \eqref{131814_19Jun24}. 
\end{Rem}
\begin{Rem}
 Thanks to Proposition \ref{3} which is stated in Section 2, Theorem \ref{2} is valid if we change $\varphi$ to any $\varphi \in {\mathcal S}(\R^n)\backslash \{0\}.$
\end{Rem}
This paper is arranged as follows: 
In Section 2, we introduce wave packet transform and its characterization of the wave front set. 
In Section 3, we briefly introduce the transformed equations of \eqref{1} via wave packet transform. 
In Section 4, we give a proof of Theorem \ref{2} which is the main result in this paper. 
In Section 5, we give a proof of Lemma 4.1 which is the key for the proof of Theorem \ref{2}. 
In Appendix A, we give the estimate of $
\left|\left(
\varphi_\lambda^{(t)},{\mathcal N}[\varphi_\lambda^{(t)}]
\right)\right|
$, which is used in Section 5. 
In Appendix B, we give some estimate for  wave packet transform of the inhomogenious term  which is \eqref{16} in Section 5. 

\section{Wave packet transform and characterization of wave front set}

In this section, we briefly give some properties of wave packet transform and introduce characterization of wave front set via wave packet transform. 
\begin{Def}[Inverse wave packet transform]
Let $F\in \mathcal{S}'(\R^n\times\R^n)$ and $\varphi\in\mathcal{S}(\R^n)\backslash\{0\}$. We define the adjoint operator $W^*_{\varphi}$ of $W_{\varphi}$ by
\[
W^*_{\varphi}\left[F(y,\xi)\right]=\frac{1}{(2\pi)^{n}}\iint_{\R^{2n}}\varphi(x-y)F(y,\xi)e^{ix\xi}dyd\xi.
\]
For $f\in\mathcal{S}'(\R^n)$ and $\varphi$, $\psi\in\mathcal{S}(\R^n)\backslash\{0\}$ satisfying $(\psi,\varphi)\neq 0$, we have inversion formula (see...
\[
f(x)=\frac{1}{(2\pi)^{n}(\psi,\varphi)}\iint_{\R^{2n}}\psi(x-y)W_{\varphi}f(y,\xi)e^{ix\xi}dyd\xi.
\]
\end{Def}
\begin{Prop}
\label{rem}
If $f$ is in $L^2(\R^n)$,
\begin{equation*}
\left\|W_{\varphi}f\right\|_{L^2(\R^{2n})}=(2\pi)^{n/2}\left\|\varphi\right\|_{L^2(\R^{n})}\left\|f\right\|_{L^2(\R^n)}
\end{equation*}
\begin{proof}[(Proof of Remark \ref{rem})]
Since $W_{\varphi}f(x,\xi)=(2\pi)^{n/2}\mathscr{F}_{y\rightarrow\xi}\left[\overline{\varphi(y-x)}f(y)\right]$, Plancherel's theorem shows
\[
\left\|W_{\varphi}f\right\|_{L^2(\R^{2n})}=(2\pi)^{n/2}\left\|\overline{\varphi(y-x)}f(y)\right\|_{L^2(\R_x^n\times\R_y^n)}=(2\pi)^{n/2}\left\|\varphi\right\|_{L^2(\R^{n})}\left\|f\right\|_{L^2(\R^n)}.
\]
\end{proof}
\end{Prop}
In \cite{Folland}, G. B. Folland has introduced the characterization of the $C^{\infty}$ wave front set in terms of wave packet transform firstly. P. G\'{e}rard \cite{Gerard} and \={O}kaji \cite{Okaji} have extended that conclusion. In \cite{chara}, K. Kato, M. Kobayashi, and S. Ito have shown the characterization without any restriction on basic wave packet.
\par
Here we introduce the characterization of the $H^s$ wave front set due K. Kato, M. Kobayashi, and S. Ito \cite{chara} in order to prove Theorem \ref{3}.
\begin{Prop}[K. Kato, M. Kobayashi, and S. Ito  ]
\label{3}
Let $s\in\R$, $0<b<1$, $(x_0,\xi_0)\in \R^n\times(\R^n\backslash\{0\})$, and $u\in\mathcal{S'}(\R^n)$. The following statements are equivalent.\\
(i)\ $(x_0,\xi_0)\notin WF_{H^s}(u)$.\\
(ii)\ There exist a neighborhood $K$ of $x_0$ and a neighborhood $V$ of $\xi_0$ such that
\[
\int_1^{\infty}\lambda^{2s+n-1}\left\|W_{{\varphi}_{\lambda}}u\left(x, \lambda\xi\right)\right\|^2_{L^2_x(K)\times L_{\xi}^2(V)}d\lambda<+\infty
\]
for all $\varphi\in\mathcal{S}(\R^n)\backslash\{0\}$, where $\varphi_{\lambda}(x)=\lambda^{nb/2}\varphi(\lambda^bx)$.\\
(iii)\ There exist $\varphi_0\in\mathcal{S}(\R^n)\backslash\{0\}$, a neighborhood $K$ of $x_0$ and a neighborhood $V$ of $\xi_0$ such that
\[
\int_1^{\infty}\lambda^{2s+n-1}\left\|W_{{\varphi}_{\lambda}}u\left(x,\lambda\xi\right)\right\|^2_{L^2_x(K)\times L_{\xi}^2(V)}d\lambda<+\infty.
\]
\end{Prop}
\begin{Rem}
 In \cite{Gerard}, P. G\'{e}rard has shown Proposition \ref{4} with basic packet $\varphi(x)=e^{-x^2}$ for $b=\frac{1}{2}$ (Proof is also in J. M. Delort \cite{Delort}). In \cite{Okaji}, \={O}kaji has shown that Proposition \ref{4}(ii) implies Proposition \ref{4}(i) if $\varphi$ satisfies $\int x^{\alpha}\varphi(x)dx\neq 0$ for some multi-index $\alpha$ for $b=\frac{1}{2}$.  Conversely \={O}kaji has shown that Proposition \ref{4}(i) implies Proposition \ref{4}(ii) if exponents of $\lambda$ are $2s+n-1-\varepsilon$ for all $\varepsilon>0$ in addition to the condition of $\varphi$. In \cite{chara}, M. Kobayashi, S. Ito and one of the authors have shown the characterization of $H^s$ wave front set without any restriction on basic wave packet.
\end{Rem}
\section{Transformed equations}
\indent
In this section, we briefly recall the transformed equation via wave packet transform which is introduced by M. Kobayashi, S. Ito and one of the authors (\cite{free}). The initial value problem (\ref{1}) is transformed to
\begin{align}
\label{4}
\begin{cases}
\left(i\partial_t+i\xi\cdot\nabla_x-\frac{1}{2}|\xi|^2\right)W_{\varphi^{(t)}}u(t,x,\xi)=W_{\varphi^{(t)}}[\mathcal N\left[ u(t)\right]](t,x,\xi),\\
W_{\varphi^{(0)}}u(0,x,\xi)=W_{\varphi}u_0(x,\xi),
\end{cases}
\end{align}
where $\varphi^{(t)}(x)=e^{it\Delta/2}[\varphi](x)$. Solving (\ref{4}) by the method of characteristics, we have the following integral equation
\begin{equation}
\label{5}
\begin{split}
&W_{\varphi^{(t)}}u(t,x,\xi)\\
=&e^{-i|\xi|^2t/2}W_{\varphi}u_0\left(x-t\xi,\xi\right)-i\int_0^t e^{-i\frac{1}{2}|\xi|^2(t-\tau)}W_{\varphi^{(\tau)}}[\mathcal N\left[u(\tau )\right]](x+(\tau-t)\xi,\xi)d\tau,
\end{split}
\end{equation}
where $x+(\tau-t)\xi$ and $\xi$ are the solutions to 
\begin{align}
\begin{cases}
\dot{x}(\tau)=\xi(\tau), \qquad x(t_0)=x,\\
\dot{\xi}(\tau)=0,\qquad \xi(t_0)=\xi.
\end{cases}
\end{align}
In the next section, we use the above integral equation \eqref{5} for the proof of Theorem \ref{2}. 
\section{Proof of Theorem \ref{2}}
\begin{proof}[(Proof of Theorem \ref{2})] Let $t_0$ be the time as in the theorem. 
We denote the following assertion by $P(\rho )$:
\newline
There exist a neighborhood $V$ of $\xi_0$ and a neighborhood $K$ of the origin 
with $K(z) = z+K$ such that 
\[
 \int_{1}^{\infty}\lambda^{2\rho+n-1}\underset{\substack{z \in\R^n\\ t\in [0, t_0]}}{\rm sup}\left\|W_{\varphi_{\lambda}^{(t-t_0)}}u\left(t,x+(t-t_0)\lambda\xi,\lambda\xi\right)\right\|^2_{L^2_x(K(z))\times L^2_{\xi}(V\cup (-V))}d\lambda<+\infty, 
\]
for $0\le t \le t_0$. 
\par
In order to prove Theorem \ref{2}, we shall show that $P(\rho)$ holds with $s\leq\rho\leq r$. 
which completes the proof with Proposition \ref{3}.  To do so, we prepare the following key lemma, which is proven in Section 5. 
\begin{Lem}
\label{4.1}
If $P(\rho)$ holds for some $\rho$ satisfying $s\leq \rho < r$, 
then we have
\begin{multline}
 \int_1^{\infty}\lambda^{2(\rho+\nu)+n-1}\\
\times\left(\int_0^{t_0}\underset{\substack{z\in\R^n}}{\rm sup}\left\|W_{\varphi_{\lambda}^{(\tau-t_0)}}[\mathcal{N}[u(\tau)]](x+(\tau-t_0)\lambda\xi,\lambda\xi)\right\|_{L^2_x(K(z))\times L^2_{\xi}(V\cup (-V))}d\tau\right)^2d\lambda
\\ <+\infty\label{145733_8Oct24}
\end{multline}
with some $\nu>0$.
\end{Lem}
Here we show the theorem assuming that Lemma \ref{4.1} is valid. 
Substituting $\varphi_{\lambda}^{(-t_0)}(x)$ for $\varphi_{\lambda}(x)$, we have
\begin{multline}
\label{6}
W_{\varphi_{\lambda}^{(t-t_0)}}u\left(t,x+(t-t_0)\lambda\xi, \lambda\xi\right)
=e^{-i|\lambda\xi|^2t/2}W_{\varphi^{(-t_0)}_{\lambda}}u_0\left(x+(t-t_0)\lambda\xi,\lambda\xi\right)\\
-i\int_0^{t_0} e^{-i\frac{1}{2}|\lambda\xi|^2(t-\tau)}W_{\varphi_{\lambda}^{(\tau-t_0)}}[\mathcal N\left[u(\tau )\right]](t,x+(\tau-t_0)\lambda\xi,\lambda\xi)d\tau .
\end{multline}
Taking $L^2(K\times (V\cup (-V)))$ norm of the both sides of the above equation \eqref{6}, we have
\begin{equation}
\label{7}
\begin{split}
&\left\|W_{\varphi_{\lambda}^{(t-t_0)}}u\left(t,x+(t-t_0)\lambda\xi, \lambda\xi\right)\right\|_{L^2_x(K)\times L^2_{\xi}(V\cup (-V))}\\
\leq&\left\|W_{\varphi^{(-t_0)}_{\lambda}}u_0\left(x+(t-t_0)\lambda\xi,\lambda\xi\right)\right\|_{L^2_x(K)\times L^2_{\xi}(V\cup (-V))}\\
&\phantom{xxxxxx}+\int_0^{t_0} \left\|W_{\varphi_{\lambda}^{(\tau-t_0)}}[\mathcal N\left[ u(\tau )\right]](x+(\tau-t_0)\lambda\xi,\lambda\xi)\right\|_{L^2_x(K)\times L^2_{\xi}(V\cup (-V))}d\tau.
\end{split}
\end{equation}
First, we show $P(-\delta)$ for any $\delta >0$ by the assumption that $u(t,\cdot ) \in H^{s,0}(\R^n)\cap H^{0,s}(\R^n)$ for $0\le t \le T$. 
Putting $\lambda\xi=\eta$, we have
\begin{equation}
\label{8}
\begin{split}
&\left\|W_{\varphi_{\lambda}^{(t-t_0)\lambda \xi}}u\left(t,x+(t-t_0),\lambda\xi\right)\right\|^2_{L^2_x(K(z))\times L^2_{\xi}(V\cup (-V))}\\
\leq&\left\|W_{\varphi_{\lambda}^{(t-t_0)}}u\left(t,x+(t-t_0)\lambda\xi,\lambda\xi\right)\right\|^2_{L^2(\mathbb{R}^{n}_{x}\times\mathbb{R}^{n}_{\xi})}\\
=&\lambda^{-n}\left\|W_{\varphi_{\lambda}^{(t-t_0)}}u\left(t,x+(t-t_0)\eta,\eta\right)\right\|^2_{L^2(\mathbb{R}^{n}_{x}\times\mathbb{R}^{n}_{\eta})} \\
=&\lambda^{-n}\left\|W_{\varphi_{\lambda}^{(t-t_0)}}u\left(t,x,\eta\right)\right\|^2_{L^2(\R^{2n})}.
\end{split}
\end{equation}
As $u(t,x)$ is in $C([-T, T];H^{s,0}(\R^n)\cap H^{0,s}(\R^n))$, Remark \ref{rem} and the conservation of $L^2$ norm of solutions to the free Schr\"{o}dinger equation yield that 
\begin{equation}
\label{10}
\begin{split}
&\lambda^{-n}\left\|W_{\varphi_{\lambda}^{(t-t_0)}}u\left(t,x,\eta\right)\right\|^2_{L^2(\R^{2n})}\\
=&\lambda^{-n}\left\|\varphi_{\lambda}^{(t-t_0)}\right\|^2_{L^2(\R^n)}\left\|u(t)\right\|^2_{L^2(\R^n)}\\
=&\lambda^{-n}\left\|\varphi\right\|^2_{L^2(\R^n)}\left\|u(t)\right\|^2_{L^2(\R^n)}.
\end{split}
\end{equation}
Hence \eqref{8} and \eqref{10} show
\begin{equation}
\label{11}
\begin{split}
&\int_{1}^{\infty}\!\!\!\lambda^{-2\delta +n-1}\underset{\substack{z\in\R^n\\ t\in [0, t_0]}}{\rm sup}\left\|W_{\varphi_{\lambda}^{(t-t_0)}}u\left(t,x+(t-t_0)\lambda\xi,\lambda\xi\right)\right\|^2_{L^2_x(K(z))\times L^2_{\xi}(V\cup (-V))}\!\!\!\!\!d\lambda\\
&\leq\left\|\varphi\right\|^2_{L^2(\R^n)}\underset{ t\in [0, t_0]}{\rm sup}\left\|u(t)\right\|^2_{L^2(\R^n)}\int_{1}^{\infty}\lambda^{-1-2\delta}d\lambda\\
&<+\infty,
\end{split}
\end{equation}
which shows $P(-\delta)$. 
\par
Assuming that $P(\rho)$ holds for $-\delta\leq \rho \leq s-\nu-n/2$, we show that $P(\rho+\nu)$. 
Let $\nu=nb(p+q+1)/2$ and a neighborhood $V'$ of $\xi_0$ satisfying
\begin{equation}
\label{12}
\int_1^{\infty}\lambda^{2(\rho+\nu)+n-1}\underset{z\in \R^n}{\mbox{sup}}\left\|W_{{\varphi}_{\lambda}^{(-t_0)}}u_0\left(x-t_0\lambda\xi, \lambda\xi)\right)\right\|^2_{L^2_x(K(z)) \times L^2_{\xi}(V'\cup (-V'))}d\lambda<+\infty.
\end{equation}
Integrating the supremum over all $z \in \R^n$ and all the time in $[0,t_0]$ of both side of the inequality \eqref{7} multiplied by $\lambda^{2(\rho+\nu)+n-1}$ from $1$ to infinity with respect to $\lambda$, we have
\begin{equation}
\label{13}
\begin{split}
&\int_1^{\infty}\lambda^{2(\rho+\nu)+n-1}\underset{\substack{z\in \R^n\\ t\in [0,t_0]}}{\mbox{sup}}\left\|W_{\varphi_{\lambda}^{(t-t_0)}}u\left(t,x+(t-t_0)\lambda\xi, \lambda\xi\right)\right\|^2_{L^2_x(K(z))\times L^2_{\xi}(V\cup (-V))}d\lambda\\
\leq&\int_1^{\infty}\lambda^{2(\rho+\nu)+n-1}\underset{z\in \R^n}{\mbox{sup}}\left\|W_{\varphi^{(-t_0)}_{\lambda}}u_0\left(x-t_0\lambda\xi,\lambda\xi\right)\right\|^2_{L^2_x(K(z))\times L^2_{\xi}(V\cup (-V))}d\lambda\\
+&\int_1^{\infty}\lambda^{2(\rho+\nu)+n-1}\left(\int_0^{t_0}\underset{\substack{z\in\R^n}}{\rm sup}\left\|W_{\varphi_{\lambda}^{(\tau-t_0)}}[\mathcal{N}[u(\tau)]](x+(\tau-t_0)\lambda\xi,\lambda\xi)\right\|_{L^2_x(K(z))\times L^2_{\xi}(V\cup (-V))}d\tau\right)^2d\lambda, 
\end{split}
\end{equation}
if $V'$ is a neighborhood of the closure of $V$. 
It follows from Lemma \ref{4.1}, \eqref{12} and \eqref{13} that $P(\rho+\nu)$.
\end{proof}
\section{Proof of Lemma \ref{4.1}}
\label{144129_18Jun24}
In the following, we write $x_{\tau}=x+(\tau-t_0)\lambda\xi$ and $t_{\tau}=\tau-t_0$ for brevity. 
\par
(Proof of Lemma \ref{4.1})
It is easy to see that for any neighborhood $\tilde{K}$ of the closure of $K$ and any neighborhood $\tilde{V}$ of the closure of $V$ there exist a positive number $C$ and a real number $\lambda_0 \ge 1$ such that
\begin{equation}
\label{16}
\begin{split}
&\left\| W_{\varphi_{\lambda}^{(t_{\tau})}}\left[{\mathcal N}[u](\tau)\right](x_\tau,\lambda\xi)\right\|_{L^2_x(K)\times L^2_{\xi}(V\cup (-V))}\\
&\le \frac{1}{\left|\left(\varphi_{\lambda}^{(t_{\tau})}, {\mathcal N}\left[\varphi_{\lambda}^{(t_{\tau})} \right]\right)\right|}
\left\| \left|W_{\varphi_{\lambda}^{(t_{\tau})}}\left[\varphi_{\lambda}^{(t_{\tau})}\right]\right|\underset{x, \xi}{\ast}\left|W_{{\mathcal N}\left[\varphi_{\lambda}^{(t_{\tau})}\right]}\left[{\mathcal N}[u](\tau)\right]\right|(x_\tau,\lambda\xi)\right\|_{L^2_x(K)\times L^2_{\xi}(V\cup (-V))}\\
\leq C&
\frac{1}{\left|\left(\varphi_{\lambda}^{(t_{\tau})}, {\mathcal N}\left[\varphi_{\lambda}^{(t_{\tau})} \right]\right)\right|}
\left\{\left\|W_{{\mathcal N}\left[\varphi_{\lambda}^{(t_{\tau})}\right]}
\left[{\mathcal N}[u](\tau)\right](x_\tau,\lambda\xi)\right\|_{L^2_x(\tilde{K})\times L^2_{\xi}(\tilde{V}\cup (-\tilde{V}))}+\lambda^{-N}\right\},
\end{split}
\end{equation}
for $\lambda \ge \lambda_0$. 
For readers' convenience, the proof is given in Appendix B.
\par
We show Lemma \ref{4.1} in the case that ${\cal N}[(u)]=u^2$ first. 
We show that $P(\rho )$ implies $P(\rho +\nu )$ with $\nu = 3nb/2$ and 
$s\le \rho \le r-\nu$. 
From the fact that
\begin{equation}
\left|\left(\varphi_\lambda^{(t)},(\varphi_\lambda^{(t)})^2
\right)
\right|\gtrsim
\begin{cases}
\lambda^{nb/2} \quad (0\le |t|\le \lambda^{-2nb})\\
\lambda^{-3nb/2} \quad (\lambda^{-2nb}\le |t|\le t_0)
\end{cases}
\label{182709_13Mar24}
\end{equation}
with $\varphi (x)= e^{-|x |^2/2}$ which is proven in Appendix A, we have
\begin{equation}
\begin{split}
&\int_0^{t_0}
\frac{1}{\left|\left(\varphi_{\lambda}^{(t_{\tau})}, \left(\varphi_{\lambda}^{(t_{\tau})} \right)^2\right)\right|}
\left\|W_{\left(\varphi_{\lambda}^{(t_{\tau})}\right)^2}\left[u^2(\tau)\right](x_\tau,\lambda\xi)\right\|_{L^2_x(\tilde{K})\times L^2_{\xi}(\tilde{V}\cup (-\tilde{V}))}d\tau\\
\lesssim&\lambda^{-nb/2}\int_{t_0-\lambda^{-2nb}}^{t_0}\left\|W_{\left(\varphi_{\lambda}^{(t_{\tau})}\right)^2}\left[u^2(\tau)\right](x_\tau,\lambda\xi)\right\|_{L^2_x(\tilde{K})\times L^2_{\xi}(\tilde{V}\cup (-\tilde{V}))}d\tau\\
&\phantom{xxxxxxxx}+\lambda^{3bn/2}\int_{0}^{t_0-\lambda^{-2nb}}\left\|W_{\left(\varphi_{\lambda}^{(t_{\tau})}\right)^2}\left[u^2(\tau)\right](x_\tau,\lambda\xi)\right\|_{L^2_x(\tilde{K})\times L^2_{\xi}(\tilde{V}\cup (-\tilde{V}))}d\tau.
\end{split}
\end{equation}

\par
Until the end of the proof in the case that ${\mathcal N}[u]=u^2$, 
we write $V$ as $\tilde{V}\cup (-\tilde{V})$ and $K$ as $\tilde{K}$. 
From the inequality \eqref{16}, it suffices to show that
\begin{equation}
\label{21}
\int_1^{\infty}\lambda^{2(\rho+\nu)+n-1+3nb}\left(\int_{0}^{t_0-\lambda^{-2nb}} \left\| W_{\left(\varphi_{\lambda}^{(t_{\tau})}\right)^2}\left[u^2(\tau)\right](x+t_{\tau}\lambda\xi, \lambda\xi) \right\|_{L^2_x(K)\times L^2_{\xi}(V)}d\tau\right)^2 d\lambda<+\infty
\end{equation}
and
\begin{equation}
\label{22}
\int_1^{\infty}\lambda^{2(\rho+\nu)+n-1-nb}\left(\int^{t_0}_{t_0-\lambda^{-2nb}} \left\| W_{\left(\varphi_{\lambda}^{(t_{\tau})}\right)^2}\left[u^2(\tau)\right](x+t_{\tau}\lambda\xi, \lambda\xi) \right\|_{L^2_x(K)\times L^2_{\xi}(V)}d\tau\right)^2 d\lambda<+\infty. 
\end{equation}
In order to show \eqref{21}, we note that 
\[
|x+t_{\tau}\lambda\xi|\geq\frac{1}{2}\lambda^{1-2nb}|\xi|
\]
for $\lambda^{-2nb}\leq|t_{\tau}|\leq t_0$, $x\in K$ and 
$\lambda\geq\lambda_0$ with some $\lambda_0$, 
since $K$ is compact and
\begin{equation*}
\begin{split}
|x+t_{\tau}\lambda\xi|&\geq|\lambda\xi(\tau-t_0)|-|x|
\geq\lambda^{1-2nb}\left(|\xi|-\frac{|x|}{\lambda^{1-2nb}}\right). 
\end{split}
\end{equation*}
Hence we have with changing variables as $x_\tau \to x$
\begin{equation}
\label{2_0717July}
\begin{split}
&\left(\int^{t_0- \lambda^{-2nb}}_{0} \left\| W_{\left(\varphi_{\lambda}^{(t_{\tau})}\right)^2}\left[u^2(\tau)\right](x_{\tau}, \lambda\xi) \right\|_{L^2_x(K)\times L^2_{\xi}(V)}d\tau\right)^2\\
=&\left(\int^{t_0- \lambda^{-2nb}}_{0} \left\|(2\pi)^{n/2}\mathscr{F}_{y\rightarrow\lambda\xi}\left[\left(\overline{\varphi_{\lambda}^{(t_{\tau})}(y-x)}u(\tau, y)\right)^2\right]\right\|_{L^2_x(K)\times L^2_{\xi}(V)}d\tau\right)^2\\
\leq&C_n\lambda^{-4s(1-2nb)}\left(\int^{t_0- \lambda^{-2nb}}_{0} \left\|\braket{x_{\tau}}^{2s}W_{\varphi_{\lambda}^{(t_{\tau})}}[u(\tau)]\underset{\xi}{\ast} W_{\varphi_{\lambda}^{(t_{\tau})}}[u(\tau)](x_{\tau}, \lambda\xi) \right\|_{L^2_x(K)\times L^2_{\xi}(V)}d\tau\right)^2 \\
\leq& C'\lambda^{-4s(1-2nb)}\left(\int^{t_0-\lambda^{-2nb}}_{0}
\left\|\left\|\braket{x}^sW_{\varphi_{\lambda}^{(t_{\tau})}}[u(\tau)](x, \xi)\right\|_{L^2(\R_{\xi}^n)}^2
\right\|_{L^4{(K)}}d\tau\right)^2\\
\leq& C'\lambda^{-4s(1-2nb)}\left(\int^{t_0-\lambda^{-2nb}}_{0}
\left\|\left\|\braket{x}^sW_{\varphi_{\lambda}^{(t_{\tau})}}[u(\tau)](x, \xi)
\right\|_{L^4{(\R^n_x)}}
\right\|_{L^2(\R_{\xi}^n)}^2
d\tau\right)^2, 
\end{split}
\end{equation}
since $V$ is bounded and does not includes the origin. Inverse wave packet transform and Young's convolution inequality yield
\begin{equation}
\label{3_0717July}
\begin{split}
&\left\|\left\|\braket{x}^sW_{\varphi_{\lambda}^{(t_{\tau})}}[u(\tau)](x, \xi)\right\|_{L^4{(\R^n_x)}}\right\|_{L^2(\R_{\xi}^n)}\\
=&\left\|\left\|\braket{x}^sW_{\varphi_{\lambda}^{(t_{\tau})}}\left[\frac{1}{(2\pi)^{n}\left(\varphi_{\lambda}^{(t_{\tau})},\varphi^{(t_{\tau})} \right)}W_{\varphi_{\lambda}^{(t_{\tau})}}^*W_{\varphi^{(t_{\tau})} }\left[u(\tau)\right]\right]\left(\tau,x,\xi\right)\right\|_{L^4{(\R^n_x)}}\right\|_{L^2(\R_{\xi}^n)}\\
\leq&\frac{C_n}{\left|\left(\varphi_{\lambda}^{(t_{\tau})},\varphi^{(t_{\tau})} \right)\right|}\left\|\left\|\left|\braket{x}^sW_{\varphi_{\lambda}^{(t_{\tau})}}\left[\varphi_{\lambda}^{(t_{\tau})}\right]\right|\underset{x, \xi}{\ast}\left|\braket{x}^sW_{\varphi^{(t_{\tau})}}\left[u(\tau)\right]\right|(x,\xi)\right\|_{L^4{(\R^n_x)}}\right\|_{L^2(\R_{\xi}^n)}\\
\leq&C_n\lambda^{nb/2}\left\|\left\|\braket{x}^sW_{\varphi_{\lambda}^{(t_{\tau})}}\left[\varphi_{\lambda}^{(t_{\tau})}\right]\right\|_{L^{4/3}(\R^n_x)}\right\|_{L^1(\R^n_{\xi})}\left\|\braket{x}^sW_{\varphi^{(t_{\tau})}}\left[u(\tau)\right]\right\|_{L^{2}(\R^{2n})}.
\end{split}
\end{equation}
The equality $W_{\varphi_{\lambda}^{(t_{\tau})}}\left[\varphi_{\lambda}^{(t_{\tau})}\right](x,\xi)=e^{-it_{\tau}|\xi|^2/2}W_{\varphi_{\lambda}}\left[\varphi_{\lambda}\right](x-t_{\tau}\xi,\xi)$ 
and the fact that 
$
W_{\varphi_{\lambda}}\left[\varphi_{\lambda}\right](x, \xi)
=
W_{\varphi}\left[\varphi\right](\lambda^bx, \lambda^{-b}\xi)
$,
yield with changing variables as $x-t_{\tau}\xi \mapsto x$
\begin{equation}
\label{4_0717July}
\begin{split}
&\left\|\left\|\braket{x}^sW_{\varphi_{\lambda}^{(t_{\tau})}}\left[\varphi_{\lambda}^{(t_{\tau})}\right]\right\|_{L^{4/3}(\R^n_x)}\right\|_{L^1(\R^n_{\xi})}\\
=&\left\|\left\|\braket{x}^se^{-it_{\tau}|\xi|^2/2}W_{\varphi_{\lambda}}\left[\varphi_{\lambda}\right](x-t_{\tau}\xi,\xi)\right\|_{L^{4/3}(\R^n_x)}\right\|_{L^1(\R^n_{\xi})}\\
=&\left\|\left\|\braket{x}^sW_{\varphi}\left[\varphi\right](\lambda^bx,\lambda^{-b}\xi)\right\|_{L^{4/3}(\R^n_x)}\right\|_{L^1(\R^n_{\xi})}\\
\leq&\lambda^{nb/4}\left\|\left\|\braket{x}^sW_{\varphi}\left[\varphi\right](x,\xi)\right\|_{L^{4/3}(\R^n_x)}\right\|_{L^1(\R^n_{\xi})}.\\
\end{split}
\end{equation}
The equalities \eqref{2_0717July}, \eqref{3_0717July} and \eqref{4_0717July} yield
\begin{equation}
\begin{split}
&\int_1^{\infty}\lambda^{2(\rho+\nu)+n-1+3nb}\left(\int^{t_0-\lambda^{-2nb}}_{0} \left\| W_{\left(\varphi_{\lambda}^{(t_{\tau})}\right)^2}\left[u^2(\tau)\right](x+t_{\tau}\lambda\xi, \lambda\xi) \right\|_{L^2_x(K)\times L^2_{\xi}(V)}d\tau\right)^2 d\lambda\\
\leq& C\int_{1}^{\infty}\lambda^{2\rho+n-1+6nb-4s(1-2bn)}\left(\int^{t_0-\lambda^{-2bn}}_{0}C\lambda^{3bn/2}d\tau\right)^2d\lambda\\
=&C \int_1^{\infty}\lambda^{2\rho+n-1-4s+bn(9+8s)}d\lambda <\infty, 
\end{split}
\end{equation}
which shows \eqref{21} if we take $b$ sufficiently small. 
\par
In order to show \eqref{22}, we devide $\R^n_\eta$ into two parts $D_3$ and $D_4$ with
$D_3=\{\eta \in \R^n | |\eta | \ge \varepsilon \lambda \text{ and } |\lambda \xi -\eta | \ge \varepsilon \lambda\}$ and 
$D_4=\{\eta \in \R^n | |\eta |\le \varepsilon \lambda \text{ or } |\lambda \xi -\eta |\le \varepsilon \lambda\}$. Since
\begin{equation}
\label{23}
\begin{split}
&W_{\left(\varphi_{\lambda}^{(t_{\tau})}\right)^2}\left[u^2(\tau)\right](x+t_{\tau}\lambda\xi, \lambda\xi)\\
=&(2\pi)^{n/2}W_{\varphi_{\lambda}^{(t_{\tau})}}[u(\tau)]\underset{\xi}{\ast} W_{\varphi_{\lambda}^{(t_{\tau})}}[u(\tau)](x+t_{\tau}\lambda\xi, \lambda\xi)\\
=& I(\tau,x,\xi,\lambda ; D_3)+I(\tau,x,\xi,\lambda ; D_4),
\end{split}
\end{equation}
where 
\begin{equation*}
\begin{split}
&I(\tau,x,\xi,\lambda ; D_j) \\
=&(2\pi)^{n/2}\int_{D_j} W_{\varphi_{\lambda}^{(t_{\tau})}}[u(\tau)](x+t_{\tau}\lambda\xi,\lambda\xi-\eta)W_{\varphi_{\lambda}^{(t_{\tau})}}[u(\tau)](x+t_{\tau}\lambda\xi,\eta)d\eta , \qquad (j = 3, 4), 
\end{split}
\end{equation*}
Schwarz's inequality and Young's inequality yield that 
\begin{align*}
\allowdisplaybreaks
\label{116}
&\left\|I(\tau,x,\xi,\lambda ; D_3)\right\|_{L^2_x(K(z))\times L^2_{\xi}(V)}\\
&\leq C \lambda^{-2s}\left\|\int_{\R^n_{\eta}}\braket{\lambda\xi-\eta}^sW_{\varphi_{\lambda}^{(t_{\tau})}}[u(\tau)]\left(x_{\tau},\lambda\xi-\eta\right)\cdot \braket{\eta}^sW_{\varphi_{\lambda}^{(t_{\tau})}}[u(\tau)]\left(x_{\tau},\eta\right)d\eta\right\|_{L^2_x(K(z))\times L^2_{\xi}(V)}\\
&\leq C \lambda^{-2s}\left\|\left\|\braket{\cdot}^sW_{\varphi_{\lambda}^{(t_{\tau})}}[u(\tau)]\left(x_{\tau},\cdot\right)\right\|^2_{L^2(\R^n)} \right\|_{L^2_x(K(z))\times L^2_{\xi}(V)}\\
&\leq C \lambda^{-2s}\left\|\left\|\braket{\xi}^sW_{\varphi_{\lambda}^{(t_{\tau})}}[u(\tau)]\left(x,\xi\right)\right\|_{L^2_\xi(\R^n)} \right\|^2_{L^4_x(\R^n)}\\
&\leq C \lambda^{-2s}\left\|\left\|\frac{1}{(2\pi)^n|(\varphi_{\lambda}^{(t_{\tau})}, \varphi^{(t_{\tau})})|}\left|\braket{\xi}^sW_{\varphi_{\lambda}^{(t_{\tau})}}[\varphi_{\lambda}^{(t_{\tau})}]\right|\underset{x, \xi}{\ast} \left|\braket{\xi}^sW_{\varphi^{(t_{\tau})}}[u(\tau)]\right|(x, \xi)\right\|_{L^2_\xi(\R^n)} \right\|^2_{L^4_x(\R^n)}\\
&\leq C\lambda^{-2s+nb}\left\|\left\| \braket{\xi}^sW_{\varphi_{\lambda}^{(t_{\tau})}}[\varphi_{\lambda}^{(t_{\tau})}] \right\|_{L^{4/3}_x(\R^n)} \right\|^2_{L^1_{\xi}(\R^n)}\left\|\left\| \braket{\xi}^sW_{\varphi^{(t_{\tau})}}[u(\tau)] \right\|_{L^{2}_x(\R^n)} \right\|^2_{L^2_{\xi}(\R^n)}\\
&= C\lambda^{-2s+nb}\left\|\left\| \braket{\xi}^sW_{\varphi}[\varphi](\lambda^bx, \lambda^{-b}\xi) \right\|_{L^{4/3}_x(\R^n)} \right\|^2_{L^1_{\xi}(\R^n)}\left\|\left\| \braket{\xi}^sW_{\varphi^{(t_{\tau})}}[u(\tau)] \right\|_{L^{2}_x(\R^n)} \right\|^2_{L^2_{\xi}(\R^n)}\\
&\leq C\lambda^{-2s+b(2s+3n/2)}\left\|\left\| \braket{\xi}^sW_{\varphi}[\varphi](x, \xi) \right\|_{L^{4/3}_x(\R^n)} \right\|^2_{L^1_{\xi}(\R^n)}\left\|\left\| \braket{\xi}^sW_{\varphi^{(t_{\tau})}}[u(\tau)] \right\|_{L^{2}_x(\R^n)} \right\|^2_{L^2_{\xi}(\R^n)}. 
\end{align*}
Hence we have
\begin{equation}
\begin{split}
&\int_1^{\infty}\lambda^{2(\rho+\nu)+n-1-nb}\Bigg(\int^{t_0}_{t_0-\lambda^{-2nb}} \left\| I(\tau,x,\xi,\lambda;D_3)\right\|_{L^2_x(K)\times L^2_{\xi}(V)}d\tau\Bigg)^2 d\lambda\\
\leq &C\int_1^{\infty}\lambda^{2\rho+n-1+2nb}\Bigg(\lambda^{-2s+b(2s-n/2)}\left\|\left\| \braket{\xi}^sW_{\varphi}[\varphi](x, \xi) \right\|_{L^{4/3}_x(\R^n)} \right\|^2_{L^1_{\xi}(\R^n)}\\
&\phantom{xxxxxxxxxxxxxxxxxxxx}\times \underset{\substack{\tau\in[0, t_0]}}{\rm sup}\left\|\left\| \braket{\xi}^sW_{\varphi^{(t_{\tau})}}[u(\tau)] \right\|_{L^{2}_x(\R^n)} \right\|^2_{L^2_{\xi}(\R^n)}\Bigg)^2d\lambda\\
=&C\int_1^{\infty}\lambda^{2\rho+n-1-4s+b(4s+n)}d\lambda <+\infty, 
\end{split}
\end{equation}
if we take $b>0$ sufficiently small as $2\rho+n-1-4s+b(4s+n)<-1$. 
\par
For the  estimate on $D_4$, it suffices to treat only on  $|\eta|\le \varepsilon \lambda $ by symmetry.
Changing variable as $\xi - \eta /\lambda $ to $\xi$ and taking a neighborhood
$V'$ of $\epsilon$-neighborhood of $V$ sufficienly small, we have
\begin{equation}
\begin{split}
\label{117}
&\left\|\int_{|\eta|\leq\varepsilon\lambda}W_{\varphi_{\lambda}^{(t_{\tau})}}[u(\tau)]\left(x_{\tau},\lambda\xi-\eta\right)W_{\varphi_{\lambda}^{(t_{\tau})}}[u(\tau)]\left(x_{\tau},\eta\right)d\eta\right\|_{L^2_x(K(z))\times L^2_{\xi}(V)}\\
\leq&\int_{|\eta|\leq\varepsilon\lambda}\left\|W_{\varphi_{\lambda}^{(t_{\tau})}}[u(\tau)]\left(x_{\tau},\lambda\xi-\eta\right)W_{\varphi_{\lambda}^{(t_{\tau})}}[u(\tau)]\left(x_{\tau},\eta\right)\right\|_{L^2_x(K(z))\times L^2_{\xi}(V)}d\eta\\
\le& \int_{|\eta|\leq\varepsilon\lambda} \Big\|W_{\varphi_{\lambda}^{(t_{\tau})}}[u(\tau)]\left(x+\eta+(\tau-t_0)\lambda\xi,\lambda\xi\right)\\
&\phantom{xxxxxxxxxxxxxx}\times W_{\varphi_{\lambda}^{(t_{\tau})}}[u(\tau)]\left(x+\eta+(\tau-t_0)\lambda\xi,\eta\right)\Big\|_{L^2_x(K(z))\times L^2_{\xi}(V')}d\eta. 
\end{split}
\end{equation}
 H\"{o}lder's inequality with respect to $x$ and $\xi$ shows that the above quanlity is estimated by
\begin{equation}
\begin{split}
 &\int_{|\eta|\leq\varepsilon\lambda}\left\|W_{\varphi_{\lambda}^{(t_{\tau})}}[u(\tau)]\left(x+\eta+(\tau-t_0)\lambda\xi,\lambda\xi\right)\right\|_{L^2_x(K(z))\times L^2_{\xi}(V')}\\
&\phantom{xxxxxxxxxxxxxx}\times \left\|W_{\varphi_{\lambda}^{(t_{\tau})}}[u(\tau)]\left(x+\eta+(\tau-t_0)\lambda\xi,\eta\right)\right\|_{L^{\infty}_x(K(z))\times L^{\infty}_{\xi}(V')} d\eta\\
\leq&\underset{\substack{z\in\R^n}}{\rm sup}\left\|W_{\varphi_{\lambda}^{(t_{\tau})}}[u(\tau)]\left(x+(\tau-t_0)\lambda\xi,\lambda\xi\right)\right\|_{L^2_x(K(z))\times L^2_{\xi}(V')}\\
&\phantom{xxxxxxxxxxxxxx}\times \int_{|\eta|\leq\varepsilon\lambda}\left\|W_{\varphi_{\lambda}^{(t_{\tau})}}[u(\tau)]\left(\cdot,\eta\right)\right\|_{L_x^{\infty}(\R^n)} d\eta\\
\leq&\underset{\substack{z\in\R^n}}{\rm sup}\left\|W_{\varphi_{\lambda}^{(t_{\tau})}}[u(\tau)]\left(x+(\tau-t_0)\lambda\xi,\lambda\xi\right)\right\|_{L^2_x(K(z))\times L^2_{\xi}(V')}\left\|W_{\varphi_{\lambda}^{(t_{\tau})}}[u(\tau)]\right\|_{L_x^{\infty}(\R^n)\times L_{\xi}^1(\R^n)}.
\end{split}
\end{equation}
Schwarz's inequality yields
\begin{equation}
\begin{split}
&\left\|W_{\varphi_{\lambda}^{(t_{\tau})}}[u(\tau)]\right\|_{L_x^{\infty}(\R^n)\times L_{\xi}^1(\R^n)}\\
=&\left\| W_{\varphi_{\lambda}^{(t_{\tau})}}\left[\frac{1}{(2\pi)^{n}\left(\varphi_{\lambda}^{(t_{\tau})}, \varphi^{(t_{\tau})}\right)}W_{\varphi_{\lambda}^{(t_{\tau})}}^*W_{ \varphi^{(t_{\tau})} }u(\tau)\right]\right\|_{L_x^{\infty}(\R^n)\times L_{\xi}^1(\R^n)}\\
\leq&\frac{C}{\left|\left(\varphi_{\lambda}^{(t_{\tau})}, \varphi^{(t_{\tau})}\right)\right|}\left\|\left|W_{\varphi_{\lambda}^{(t_{\tau})}}[\varphi_{\lambda}^{(t_{\tau})}]\right|\underset{x,\xi}{\ast}\left|W_{\varphi{(t_{\tau})}}[u(\tau)]\right|\right\|_{L_x^{\infty}(\R^n)\times L_{\xi}^1(\R^n)}\\
\leq & C\lambda^{nb/2}\left\|W_{\varphi_{\lambda}^{(t_{\tau})}}[\varphi_{\lambda}^{(t_{\tau})}]\right\|_{L_x^2(\R^n)\times L_{\xi}^1(\R^n)}\left\|W_{\varphi{(t_{\tau})}}[u(\tau)]\right\|_{L_x^2(\R^n)\times L_{\xi}^1(\R^n)}\\
\leq  & C\lambda^{nb/2}\left\|W_{\varphi_{\lambda}^{(t_{\tau})}}[\varphi_{\lambda}^{(t_{\tau})}]\right\|_{L_x^2(\R^n)\times L_{\xi}^1(\R^n)}\left\|\braket{\xi}^{-s}\right\|_{L^2(\R^n)}\left\|\braket{\xi}^sW_{\varphi{(t_{\tau})}}[u(\tau)]\right\|_{L_x^2(\R^n)\times L_{\xi}^2(\R^n)}
\end{split}
\end{equation}
The equality $W_{\varphi_{\lambda}^{(t_{\tau})}}\left[\varphi_{\lambda}^{(t_{\tau})}\right](x,\xi)=e^{-it_{\tau}|\xi|^2/2}W_{\varphi_{\lambda}}\left[\varphi_{\lambda}\right](x-t_{\tau}\xi,\xi)$ 
and the fact that 
$
W_{\varphi_{\lambda}}\left[\varphi_{\lambda}\right](x, \xi)
=
W_{\varphi}\left[\varphi\right](\lambda^bx, \lambda^{-b}\xi)
$,
yield with changing variables as $x-t_{\tau}\xi \mapsto x$
\begin{equation}
\begin{split}
\left\|W_{\varphi_{\lambda}^{(t_{\tau})}}[\varphi_{\lambda}^{(t_{\tau})}]\right\|_{L_x^2(\R^n)\times L_{\xi}^1(\R^n)}&=\left\|W_{\varphi_{\lambda}}[\varphi_{\lambda}](x-t_{\tau}\xi, \xi)\right\|_{L_x^2(\R^n)\times L_{\xi}^1(\R^n)}\\
&=\left\|W_{\varphi}[\varphi](\lambda^bx, \lambda^{-b}\xi)\right\|_{L_x^2(\R^n)\times L_{\xi}^1(\R^n)}\\
&=\lambda^{nb/2}\left\|W_{\varphi}[\varphi]\right\|_{L_x^2(\R^n)\times L_{\xi}^1(\R^n)}.\\
\end{split}
\end{equation}
Hence we have
\begin{equation}
\label{072401}
\begin{split}
&\int_1^{\infty}\lambda^{2(\rho+\nu)+n-1-nb}\Bigg(\int_{t_0-\lambda^{-2bn}}^{t_0}\Bigg\|\int_{|\eta|\leq\varepsilon\lambda}W_{\varphi_{\lambda}^{(t_{\tau})}}[u(\tau)]\left(x_{\tau},\lambda\xi-\eta\right)\\
&\phantom{xxxxxxxxxxxxxxxxxxxxxxxxxxxxxx}\times W_{\varphi_{\lambda}^{(t_{\tau})}}[u(\tau)]\left(x_{\tau},\eta\right)d\eta\Bigg\|_{L^2_x(K(z))\times L^2_{\xi}(V)}d\tau\Bigg)^2d\lambda\\
\leq&\int_1^{\infty}\lambda^{2\rho+n-1+2nb}\left(C\lambda^{nb}\int_{t_0-\lambda^{-2nb}}^{t_0}\underset{\substack{z\in\R^n}}{\rm sup}\left\|W_{\varphi_{\lambda}^{(t_{\tau})}}[u(\tau)]\left(x+(\tau-t_0)\lambda\xi,\lambda\xi\right)\right\|_{L^2_x(K(z))\times L^2_{\xi}(V')}\right)^2d\lambda\\
\leq&\int_1^{\infty}\lambda^{2\rho+n-1}\underset{\substack{z\in\R^n\\ t\in[0, t_0]}}{\rm sup}\left\|W_{\varphi_{\lambda}^{(t_{\tau})}}[u(\tau)]\left(x+(\tau-t_0)\lambda\xi,\lambda\xi\right)\right\|_{L^2_x(K(z))\times L^2_{\xi}(V')}d\tau\\
<&+\infty , 
\end{split}
\end{equation}
if we take $V$ sufficiently small as $P(\rho )$ holds. 
On $|\lambda\xi-\eta|\leq \varepsilon\lambda$, the inequality \eqref{072401} also holds  with changing variables as $\lambda\xi-\eta \mapsto \eta$. Thus \eqref{22} holds.
\par
Since 
\begin{equation}
 W_{\varphi^{(t)}_\lambda}[\overline{u}](x,\xi ) 
= \overline{
W_{\varphi^{(-t)}_\lambda}u(x,-\xi )
}, 
\end{equation} 
the proof in the case that ${\mathcal N}[u]=\overline{u}u$ 
can be done by the same way as in the case that ${\mathcal N}[u]=u^2$ 
with the help of the condition \eqref{1_July17} and the fact that 
${\mathcal N}[\varphi^{(t)}] = \varphi^{(t)}\overline{\varphi^{(t)}}= \varphi^{(t)}\varphi^{(-t)} $. 
The proof in the case that ${\mathcal N}[u]=u^p\overline{u}^q$ can
be completed by induction. \qed
%
%
\section*{Appendix A}
We give a proof of the fact that
\begin{equation}
 \left|\left(\varphi_\lambda^{(t)},{\mathcal N}[\varphi_\lambda^{(t)}]
\right)
\right|\gtrsim
\begin{cases}
\lambda^{nb(p+q-1)/2} \quad (0\le |t|\le \lambda^{-(p+q)nb})\\
\lambda^{-nb(p+q+1)/2} \quad (\lambda^{-(p+q)nb}\le |t|\le t_0)
\end{cases}
\label{090401}
\end{equation}
with $\varphi (x)= e^{-|\xi |^2/2}$. 
\begin{proof}
From the fact that 
$$
\varphi_\lambda^{(t)}= c_n \lambda^{nb/2}(1+\lambda^{4b}t^2) ^{-n/4}e^{-\frac{(1-i\lambda^{2b}t)}{2(1+\lambda^{4b}t^2)}|\lambda^bx|^2}, 
$$
we have 
\begin{align*}
\left|\left(\varphi_\lambda^{(t)},{\mathcal N}[\varphi_\lambda^{(t)}]\right)\right|
&=
\left|
\lambda^{nb(p+q-1)/2}(1+\lambda^{4b}t^2)^{-n(p+q+1)/4}
\int e^{-\frac{(p+q-1)+i(p-q-1)\lambda^{2b}t}{2(1+\lambda^{4b}t^2)}|x|^2}\, dx
\right|\\
& = 
\lambda^{nb(p+q-1)/2}(1+\lambda^{4b}t^2)^{-n(p+q+1)/4}
\left|
\frac{2(1+\lambda^{4b}t^2)}{(p+q-1)+i(p-q-1) \lambda^{2b}t}\pi
\right|^{n/2}
\\
&\ge C \lambda^{nb(p+q-1)/2}(1+\lambda^{4b}t^2)^{-n(p+q)/4}, 
\end{align*}
which yields \eqref{182709_13Mar24}. 
\end{proof}
We give a proof of the fact that
\begin{equation}
 \left|\left(\varphi_\lambda^{(t)},(\varphi_\lambda^{(t)})^2
\right)
\right|\gtrsim
\begin{cases}
\lambda^{nb/2} \quad (0\le |t|\le \lambda^{-2nb})\\
\lambda^{-3nb/2} \quad (\lambda^{-2nb}\le |t|\le t_0)
\end{cases}
\label{182709_13Mar24}
\end{equation}
with $\varphi (x)= e^{-|\xi |^2/2}$. 
\begin{proof}
From the fact that 
$$
\varphi_\lambda^{(t)}= c_n \lambda^{nb/2}(1+\lambda^{4b}t^2) ^{-n/4}e^{-\frac{(1-i\lambda^{2b}t)}{2(1+\lambda^{4b}t^2)}|x|^2}, 
$$
we have 
\begin{align*}
|(\varphi_\lambda^{(t)},(\varphi_\lambda^{(t)})^2)|
&=
\left|
\lambda^{nb/2}(1+\lambda^{4b}t^2)^{-3n/4}
\int e^{-\frac{3+i\lambda^{2b}t}{2(1+\lambda^{4b}t^2)}|x|^2}\, dx
\right|\\
& = 
\lambda^{nb/2}(1+\lambda^{4b}t^2)^{-3n/4}
\left|
\frac{2(1+\lambda^{4b}t^2)}{3+i \lambda^{2b}t}\pi
\right|^{n/2}
\\
&\ge C \lambda^{nb/2}(1+\lambda^{4b}t^2)^{-n/2}, 
\end{align*}
which yields \eqref{182709_13Mar24}. 
\end{proof}
\section*{Appendix B}
We show \eqref{16} with the notation that $f(t,x)={\mathcal N}[u] \in C([0,T], H^{s,0}\cap H^{0,s})$. 
\begin{equation}
\label{14}
\begin{split}
&\left\|W_{\varphi_{\lambda}^{(\tau-t_0)}}[f(\tau)](x_\tau,\lambda\xi)\right\|_{L^2_x(K)\times L^2_{\xi}(V)}\\
=&\left\|W_{\varphi_{\lambda}^{(t_{\tau})}}\left[\frac{1}{(2\pi)^{n}\left(\varphi_{\lambda}^{(t_{\tau})}, {\mathcal N}[\varphi_{\lambda}^{(t_\tau)}]\right)}W_{\varphi_{\lambda}^{(t_{\tau})}}^*W_{{\mathcal N}[\varphi_{\lambda}^{(t_\tau)}]}\left[f(\tau)\right]\right]\left(\tau,x_{\tau},\lambda\xi\right)\right\|_{L^2_x(K)\times L^2_{\xi}(V)}\\
\leq&\displaystyle C\frac{1}{\left|\left(\varphi_{\lambda}^{(t_{\tau})}, {\mathcal N}[\varphi_{\lambda}^{(t_\tau)}]\right)\right|}\Bigg\| \left|W_{\varphi_{\lambda}^{(t_{\tau})}}\left[\varphi_{\lambda}^{(t_{\tau})}\right]\right|\underset{x, \xi}{\ast}\left|W_{\left(\varphi_{\lambda}^{(t_{\tau})}\right)^2}\left[f(\tau)\right]\right|(x_\tau,\lambda\xi)\Bigg\|_{L^2_x(K)\times L^2_{\xi}(V)}.\\ 
\end{split}
\end{equation}
The equality $W_{\varphi_{\lambda}^{(t_{\tau})}}\left[\varphi_{\lambda}^{(t_{\tau})}\right](x,\xi)=e^{-it_{\tau}|\xi|^2/2}W_{\varphi_{\lambda}}\left[\varphi_{\lambda}\right](x-t_{\tau}\xi,\xi)$ 
and the fact that 
$
W_{\varphi_{\lambda}}\left[\varphi_{\lambda}\right](w, \eta)
=
W_{\varphi}\left[\varphi\right](\lambda^bw, \lambda^{-b}\eta)
$,
yield with changing variables as $w-t_{\tau}\eta \mapsto w$
\begin{equation}
\label{18}
\begin{split}
& \left|W_{\varphi_{\lambda}^{(t_{\tau})}}\left[\varphi_{\lambda}^{(t_{\tau})}\right]\right|\underset{x, \xi}{\ast}\left|W_{{\mathcal N}[\varphi_{\lambda}^{(t_\tau)}]}\left[f(\tau)\right]\right|(x_\tau,\lambda\xi) \\
=&\iint\left|W_{\varphi_{\lambda}}\left[\varphi_{\lambda}\right](w-t_{\tau}\eta, \eta)\right|\cdot\left| W_{{\mathcal N}[\varphi_{\lambda}^{(t_\tau)}]}\left[f(\tau)\right](x_{\tau}-w, \lambda\xi-\eta) \right|dwd\eta \\
=&\iint\left|W_{\varphi}\left[\varphi\right](\lambda^b w, \lambda^{-b}\eta)\right|\cdot\left| W_{{\mathcal N}[\varphi_{\lambda}^{(t_\tau)}]}\left[f(\tau)\right](x-w+t_\tau (\lambda\xi-\eta), \lambda\xi-\eta) \right|dwd\eta. 
\end{split}
\end{equation}
\par
We devide ${\mathbb R}^{2n}$ into two parts 
$D_1 =\{ (x,w)\in{\mathbb R}^{2n}| |w|\ge \lambda^{-b/2} \text{ or } 
|\eta| \ge \lambda^{3b/2}\}$ 
and $D_2=\{ (x,w)\in{\mathbb R}^{2n}| |w|\le \lambda^{-b/2} \text{ and } 
|\eta| \le \lambda^{3b/2}\}$. 
Since $\braket{\lambda^b w;\lambda^{-b}\eta}\ge \lambda^{b/2}$ in $D_1$, we have
\begin{equation}
\label{19}
\begin{split}
&\iint_{D_1}\left|W_{\varphi}\left[\varphi\right](\lambda^bw, \lambda^{-b}\eta)\right|
 \left| W_{{\mathcal N}[\varphi_{\lambda}^{(t_\tau)}]}\left[f(\tau)\right](x-w+t_{\tau}(\lambda \xi -\eta), \lambda\xi-\eta) \right|dwd\eta\\  
\leq&C\lambda^{-bN/2}\int_{\R^n_{\eta}}\int_{\R^n_{w}}\braket{\lambda^b w;\lambda^{-b}\eta}^N\left|W_{\varphi}\left[\varphi\right](\lambda^bw, \lambda^{-b}\eta)\right|\\
&\phantom{xxxxxxxxxxx}\times \left| W_{{\mathcal N}[\varphi_{\lambda}^{(t_\tau)}]}\left[f(\tau)\right](x-w+t_{\tau}(\lambda \xi -\eta) , \lambda\xi-\eta) \right|dwd\eta\\
\leq&C\lambda^{-bN/2}\left\|\braket{\cdot;\cdot}^N W_{\varphi}[\varphi]\right\|_{L^2(\R^{2n})}\left\|W_{{\mathcal N}[\varphi_{\lambda}^{(t_\tau)}]}\left[f(\tau)\right]\right\|_{L^2(\R^{2n})}\\
\le& C\lambda^{-bN/2}\left\|\braket{\cdot;\cdot}^N W_{\varphi}[\varphi]\right\|_{L^2(\R^{2n})}
\left\|{\mathcal N}[\varphi_{\lambda}^{(t_\tau)}] \right\|_{L^2(\R^n)}
\left\|f\right\|_{L^2(\R^{n})}\\
\leq&C\lambda^{-bN/2+nb(p+q-1)/2},
\end{split}
\end{equation}
since $\left\|{\mathcal N}[\varphi_{\lambda}^{(t_{\tau})}] \right\|_{L^{2}(\R^n)}$ is estimated by
$$\left\|{\mathcal F}{\mathcal N}[\varphi_{\lambda}^{(t_{\tau})}] \right\|_{L^{2}(\R^n)}=
\left\|
{\mathcal F}[\varphi_{\lambda}^{(t_{\tau})}]* \cdots * {\mathcal F}[\overline{\varphi_{\lambda}^{(t_{\tau})}}]
\right\|_{L^2}
\le 
\left\|
\hat{\varphi}_\lambda
\right\|_{L^1}^{p+q-1}
\left\|\hat{\varphi}_{\lambda} \right\|_{L^2}
\le C\lambda^{nb(p+q-1)/2}. 
$$
Changing variables as $\xi-\frac{\eta}{\lambda}\mapsto \xi$, we have
\begin{equation}
\label{20}
\begin{split}
&\Bigg\|\iint_{D_2}\left|W_{\varphi}\left[\varphi\right](\lambda^b w, \lambda^{-b}\eta)\right|\\
&\phantom{xxxxxx}\times \left| W_{{\mathcal N}[\varphi_{\lambda}^{(t_\tau)}]}\left[f(\tau)\right](x-w+t_\tau (\lambda\xi-\eta), \lambda\xi-\eta) \right|dwd\eta\Bigg\|_{L^2_x(K)\times L^2_{\xi}(V)}\\
\leq& \iint\left|W_{\varphi}\left[\varphi\right](\lambda^b w, \lambda^{-b}\eta)\right|\\
&\phantom{xxxxxx}\times\left\| W_{{\mathcal N}[\varphi_{\lambda}^{(t_\tau)}]}\left[f(\tau)\right](x-w+t_{\tau}(\lambda\xi-\eta), \lambda\xi-\eta) \right\|_{L^2_x(K)\times L^2_{\xi}(V)}dwd\eta\\
\leq&\iint\left|W_{\varphi}\left[\varphi\right](\lambda^b w, \lambda^{-b}\eta)\right|\cdot\left\| W_{{\mathcal N}[\varphi_{\lambda}^{(t_\tau)}]}\left[f(\tau)\right](x+t_{\tau}\lambda\xi, \lambda\xi) \right\|_{L^2_x(\tilde{K})\times L^2_{\xi}(\tilde{V})}dwd\eta\\
\leq& \left\| W_{{\mathcal N}[\varphi_{\lambda}^{(t_\tau)}]}\left[f(\tau)\right](x+t_{\tau}\lambda\xi, \lambda\xi) \right\|_{L^2_x(\tilde{K})\times L^2_{\xi}(\tilde{V})} \left\|W_{\varphi}\left[\varphi\right]\right\|_{L^1(\R^n)},
\end{split}
\end{equation}
which shows with \eqref{19}
\begin{equation}
\label{eq-16}
\begin{split}
&\left\| \left|W_{\varphi_{\lambda}^{(t_{\tau})}}\left[\varphi_{\lambda}^{(t_{\tau})}\right]\right|\underset{x, \xi}{\ast}\left|W_{{\mathcal N}[\varphi_{\lambda}^{(t_\tau)}]}\left[f(\tau)\right]\right|(x_\tau,\lambda\xi)\right\|_{L^2_x(K)\times L^2_{\xi}(V)}\\
\leq C&\left\{\left\|W_{{\mathcal N}[\varphi_{\lambda}^{(t_\tau)}}]
 \left[f(\tau)\right](x_\tau,\lambda\xi)\right\|_{L^2_x(\tilde{K})\times L^2_{\xi}(\tilde{V})}+\lambda^{-N}\right\}.
\end{split}
\end{equation}

\end{document}